\documentstyle[12pt, amsfonts]{amsart}
\begin{document}
\def\R{{\mathbb R}}
\def\Z{{\mathbb Z}}
\def\C{{\mathbb C}}
\newcommand{\trace}{\rm trace}
\newcommand{\Ex}{{\mathbb{E}}}
\newcommand{\Prob}{{\mathbb{P}}}
\newcommand{\E}{{\cal E}}
\newcommand{\F}{{\cal F}}
\newtheorem{df}{Definition}
\newtheorem{theorem}{Theorem}
\newtheorem{lemma}{Lemma}
\newtheorem{pr}{Proposition}
\newtheorem{co}{Corollary}
\def\n{\nu}
\def\sign{\mbox{ sign }}
\def\a{\alpha}
\def\N{{\mathbb N}}
\def\A{{\cal A}}
\def\L{{\cal L}}
\def\X{{\cal X}}
\def\F{{\cal F}}
\def\c{\bar{c}}
\def\v{\nu}
\def\d{\delta}
\def\diam{\mbox{\rm dim}}
\def\vol{\mbox{\rm Vol}}  
\def\b{\beta}
\def\t{\theta}
\def\l{\lambda}
\def\e{\varepsilon}
\def\colon{{:}\;}
\def\pf{\noindent {\bf Proof :  \  }}
\def\endpf{ \begin{flushright}
$ \Box $ \\
\end{flushright}}
%%%%%%%%%%%%%%%%%%%%%%%%%%%%%%%%%%%%%%%%%%%%%%%%%%%%%%%%%%%%%%%%%%%%

\title[Complex Busemann-Petty problem]{The complex Busemann-Petty
problem on sections of convex bodies}

\author{A. Koldobsky}

\address{
Alexander Koldobsky\\
Department of Mathematics\\
University of Missouri\\
Columbia, MO 65211, USA}

\email{koldobsk@@math.missouri.edu} 

\author{H. K\"onig}
\address{
Hermann K\"onig\\
Mathematisches Seminar \\
Christian-Albrechts-UniversitŠt Kiel\\
Ludewig-Meyn Str. 4\\
D-24098 Kiel}

\email{hkoenig@@math.uni-kiel.de}

\author{M. Zymonopoulou}
\address{
Marisa Zymonopoulou\\
Department of Mathematics\\
University of Missouri\\
Columbia, MO 65211, USA}

\email{marisa@@math.missouri.edu}
%%%%%%%%%%%%%%%%%%%%%%%%%%%%%%%%%%%%%%%%%%%%%%%%%%%%%%%%%%%%%%%%%%%%
\begin{abstract} The complex Busemann-Petty problem asks whether origin 
symmetric convex bodies in $\C^n$ with smaller central hyperplane sections 
necessarily have smaller volume. We prove that the answer
is affirmative if $n\le 3$ and negative if $n\ge 4.$
\end{abstract}
\maketitle
%%%%%%%%%%%%%%%%%%%%%%%%%%%%%%%%%%%%%%%%%%%%%%%%%%%%%%%%%%%%%%%%%%%%

\section{Introduction}

The Busemann-Petty problem, posed in 1956 (see \cite{BP}), 
asks the following question. Suppose that $K$ and $L$ are origin symmetric 
convex bodies in $\R^n$ such that 
$$\vol_{n-1}(K\cap H) \le \vol_{n-1}(L\cap H)$$
for every central hyperplane $H$ in $\R^n.$ Does it follow that 
$$\vol_n(K) \le \vol_n(L)?$$
The answer is affirmative if $n\le 4$ and negative if $n\ge 5.$
The solution was completed in the end of the 90's as the result 
of a sequence of papers \cite{LR}, \cite{Ba}, \cite{Gi}, 
\cite{Bo}, \cite{Lu}, \cite{Pa}, \cite{Ga1}, \cite{Ga2}, \cite{Zh1}, 
\cite{K1}, \cite{K2}, \cite{Zh2}, \cite{GKS} ; see \cite[p. 3]{K-book} for the 
history of the solution.

In this article we consider the complex version of the problem.  For $\xi\in \C^n,
|\xi|=1,$ denote by 
$$H_\xi = \{ z\in \C^n:\ (z,\xi)=\sum_{k=1}^n z_k\overline{\xi_k} =0\}$$
the complex hyperplane perpendicular to $\xi.$

Origin symmetric convex bodies in $\C^n$ are the unit balls of norms on $\C^n.$ We denote by $\|\cdot\|_K$
the norm corresponding to the body $K:$
$$K=\{z\in \C^n:\ \|z\|_K\le 1\}.$$
In order to define volume, we identify $\C^n$ with $\R^{2n}$ using the mapping
$$\xi = (\xi_1,...,\xi_n)=(\xi_{11}+i\xi_{12},...,\xi_{n1}+i\xi_{n2})
 \mapsto  (\xi_{11},\xi_{12},...,\xi_{n1},\xi_{n2}).$$
 Under this mapping the hyperplane $H_\xi$ turns into a two-codimensional subspace of
 $\R^{2n}$ orthogonal to the vectors 
 $$\xi=(\xi_{11},\xi_{12},...,\xi_{n1},\xi_{n2})\quad  {\rm and } \quad 
 \xi^\perp =(-\xi_{12},\xi_{11},...,-\xi_{n2},\xi_{n1}).$$
 Since norms on $\C^n$ satisfy the equality 
$$\|\lambda z\| = |\lambda|\|z\|,\quad \forall z\in \C^n,\  \forall\lambda \in \C,$$
origin symmetric complex convex bodies correspond to those origin symmetric convex bodies 
$K$  in $\R^{2n}$ that are invariant 
 with respect to any coordinate-wise two-dimensional rotation, namely for each $\theta\in [0,2\pi]$
 and each $\xi= (\xi_{11},\xi_{12},...,\xi_{n1},\xi_{n2})\in \R^{2n}$
  \begin{equation} \label{rotation}
  \|\xi\|_K = 
 \|R_\theta(\xi_{11},\xi_{12}),...,R_\theta(\xi_{n1},\xi_{n2})\|_K,
 \end{equation}
 where $R_\theta$ stands for  the counterclockwise rotation of $\R^2$ by the angle 
 $\theta$ with respect to the origin. We shall simply say that $K$ {\it is invariant with respect 
 to all $R_\theta$} if it satisfies the equations  (\ref{rotation}).
 
 Now the complex Busemann-Petty problem can be formulated as follows: 
 suppose $K$ and $L$ are origin symmetric invariant with respect to all $R_\theta$ 
 convex bodies in $\R^{2n}$ such that
 $$\vol_{2n-2}(K\cap H_\xi)\le \vol_{2n-2}(L\cap H_\xi)$$
 for each $\xi$ from the unit sphere $S^{2n-1}$ of $\R^{2n}.$ Does it follow that
 $$\vol_{2n}(K) \le \vol_{2n}(L) ?$$
 
 This formulation reminds of the lower-dimensional Busemann-Petty problem, where 
 one tries to deduce the inequality for $2n$-dimensional volumes of arbitrary origin-symmetric 
 convex bodies from the inequalities for volumes of all $(2n-2)$-dimensional sections. In the case 
 where $n=2$ this amounts to  considering  two-dimensional sections of four-dimensional bodies, 
 where the answer to the lower dimensional problem is
 affirmative by the solution to the original Busemann-Petty problem - we first get
 inequalities for the volumes of all three-dimensional sections and then the inequality 
 for the four-dimensional volumes. However, if $n=3$ we get four-dimensional sections 
 of six-dimensional bodies, where the answer to the lower-dimensional 
 problem is negative by a result of Bourgain and Zhang \cite{BZ}.  Our problem is
 different from the lower-dimensional Busemann-Petty problem in two aspects.
 First, we do not have all $(2n-2)$-dimensional sections, we only have sections 
 by subspaces coming from complex hyperplanes, which makes the situation worse than
 for the lower-dimensional problem. Secondly, we consider only those convex bodies in $\R^{2n}$
 that are invariant with respect to all $R_\theta$, and we may be able to convert this invariance
 into affirmative answer in some higher dimensions. 
 
 The latter appears to be the case, as we prove below that the answer 
 to the complex Busemann-Petty problem is affirmative if $n\le 3$ and negative if $n\ge 4.$
 
 In 1988 Lutwak \cite{Lu} introduced the class of intersection bodies and found a connection 
 between this class and the ``real" Busemann-Petty problem, which played an important role in the
 solution of the problem. It appears that the complex Busemann-Petty problem is closely
 related to the class of 2-intersection bodies introduced in \cite{K5, K8}, namely the answer to the problem 
 is affirmative if and only if every origin symmetric invariant with respect to all $R_\theta$ convex body 
 in $\R^{2n}$ is a $2$-intersection body.
 We shall prove this connection in Theorem \ref{connection}.  After that we prove that 
 every origin symmetric invariant with respect to all $R_\theta$ convex body in $\R^{2n}$ 
 is a $(2n-4)$-intersection body, but not every such body is a $(2n-6)$-intersection body. 
 Putting $n=3$ and then $n=4$, one can see how these results imply the solution of the complex 
 Busemann-Petty problem. Our proofs use several results from the recently developed  Fourier
 analytic approach to sections of convex bodies; see \cite{K-book}. In Section 2, we collect 
 necessary definitions and results related to this approach.
 
 For other results related to the Busemann-Petty problem see \cite{BZ}, \cite{BFM}, \cite{K5}, \cite{K9},  \cite{KYY},  
 \cite{Mi}, \cite{Ru}, \cite{RZ}, \cite{Y1},  \cite{Y2}, \cite{Zv1}, \cite{Zv2}.
 
  \section{Elements of the Fourier approach to sections}
 
 Our main tool is the Fourier transform of distributions. As usual, we denote by ${\cal{S}}(\R^n)$ 
the Schwartz space of rapidly decreasing infinitely differentiable 
functions (test functions) in $\R^n,$ and 
${\cal{S}}^{'}(\R^n)$ is
the space of distributions over ${\cal{S}}(\R^n).$  
The Fourier transform $\hat{f}$ of a distribution $f\in {\cal{S}}^{'}(\R^n)$ 
is defined by $\langle \hat{f},\phi \rangle = 
\langle f,\hat\phi \rangle$ for every
test function $\phi.$  A distribution is called even homogeneous
of degree $p\in \R$  if 
$\langle f(x), \phi(x/\alpha) \rangle = |\alpha|^{n+p} 
\langle f,\phi \rangle$
for every test function $\phi$ and every 
$\alpha\in \R,\ \alpha\neq 0.$ The Fourier transform of 
an even homogeneous distribution of degree $p$ is an even 
homogeneous distribution of degree $-n-p.$  A distribution $f$ is called 
positive definite if, for every test function $\phi,$ 
$\langle f, \phi * \overline{\phi(-x)} \rangle \ge 0.$
This is equivalent to $\hat{f}$ being a positive distribution
in the sense that $\langle \hat{f}, \phi \rangle \ge 0$ for
every non-negative test function $\phi.$

A compact set $K$ in $\R^n$ is called a {\it star body} if every straight
 line through the origin crosses the boundary at exactly two points different from the origin, and 
 the boundary of $K$ is continuous in the sense that the {\it Minkowski functional} of $K$ defined
 by 
 $$\|x\|_K= \min\{a\ge 0:\ x\in aK\}$$
 is a continuous function on $\R^n.$ If in addition $K$ is origin symmetric and convex, then
 the Minkowski functional is a norm on $\R^n.$ If $\xi\in S^{n-1},$ then $\rho_K(\xi)=\|\xi\|_K^{-1}$
 is the radius of $K$ in the direction $\xi.$

 A simple calculation in polar coordinates gives the following {\it polar formula for 
 the volume}:
 $$n\ \vol_n(K) = n \int_{\R^n} \chi(\|x\|_K)\ dx = \int_{S^{n-1}} \|\xi\|_K^{-n}\ d\xi,$$
 where $\chi$ is the indicator function of the interval $[0,1].$

We say that a star body $K$ in $\R^n$  is $k$-smooth (infinitely smooth)
 if the restriction of $\|x\|_K$ to the sphere $S^{n-1}$ belongs to the class $C^k(S^{n-1})$ ($C^\infty(S^{n-1})$) of
 $k$ times continuously differentiable (infinitely differentiable) functions on the sphere. 
 It is well-known that one can approximate any
 convex body in $\R^n$ in the radial metric
 $$d(K,L) = \sup_{\xi\in S^{n-1}}|\rho_K(\xi)-\rho_L(\xi)|$$
 by a sequence of infinitely smooth convex bodies. This can be proved by a simple convolution
 argument (see for example \cite[Th. 3.3.1]{Sch}). It is also easy to see that any convex body in
 $\R^{2n}$ invariant with respect to all $R_\theta$ can be approximated in the radial metric 
 by a sequence of infinitely smooth convex bodies invariant with respect to all $R_\theta.$
 This follows from the same convolution argument, because invariance 
 with respect to $R_\theta$ is preserved under convolutions.

 As proved in \cite[Lemma 3.16]{K-book},  if 
 $K$ is an infinitely smooth origin symmetric star body in $\R^n$  and $0<p<n$ then
 the Fourier transform of the distribution $\|x\|_K^{-p}$ is a homogeneous function of degree 
 $-n+p$ on $\R^n,$ whose restriction to the sphere  is infinitely smooth.
 We use a version of Parseval's formula on the sphere established in \cite{K5} (see also 
 \cite[Lemma 3.22]{K-book}):
 \begin{pr} \label{parseval} Let $K$ and $L$ be infinitely smooth origin symmetric star
 bodies in $\R^n$ and $0<p<n.$ Then
 $$\int_{S^{n-1}} \left(\|x\|_K^{-p}\right)^\wedge(\xi) \left(\|x\|_L^{-n+p}\right)^\wedge(\xi)\ d\xi
 = (2\pi)^n \int_{S^{n-1}} \|x\|_K^{-p} \|x\|_L^{-n+p}\ dx.$$
  \end{pr}
 \bigbreak
The classes of $k$-intersection bodies were introduced in \cite{K5}, \cite{K8} as follows.
Let $1\le k < n,$ and let $D$ and $L$ be origin symmetric star 
bodies in $\R^n.$ We say that $D$ is 
a {\it $k$-intersection body of}  $L$ if for every $(n-k)$-dimensional subspace $H$ of $\R^n$
$$\vol_k(D\cap H^\bot) = \vol_{n-k}(L\cap H).$$
More generally, we say that an origin symmetric star body $D$ in $\R^n$ is a {\it $k$-intersection body} 
if there exists a finite Borel measure $\mu$ on $S^{n-1}$ so that for every even test function $\phi\in {\cal{S}}(\R^n),$ 
$$\int_{\R^n} \|x\|_D^{-k} \phi(x)\ dx = 
\int_{S^{n-1}} \left(
\int_0^\infty t^{k-1} \hat\phi(t\xi)\ dt\right) d\mu(\xi).$$
Note that $k$-intersection bodies of star bodies are those $k$-intersection
bodies for which the measure $\mu$ has a continuous strictly positive density; see \cite{K8} or \cite[p. 77]{K-book}.
When $k=1$ we get the class of intersection bodies introduced by Lutwak in \cite{Lu}.

A more general concept of embedding in $L_{-p}$ was introduced in \cite{K7}. 
Let $D$ be an origin symmetric star body in $\R^n,$ and 
$X=(\R^n,\|\cdot\|_D).$ For $0<p<n,$ 
we say that $X$ embeds in $L_{-p}$
if there exists a finite Borel measure $\mu $ on $S^{n-1}$ so that, for every even 
test function $\phi$
$$
\int_{\R^n} \|x\|_D^{-p} \phi(x)\ dx = \int_{S^{n-1}} \left( 
\int_{\R} |z|^{p-1} \hat{\phi}(z\theta)\ dz\right) d\mu(\theta).
$$
Obviously, an origin symmetric star body $D$ in $\R^n$ is a $k$-intersection body
if and only if the space $(\R^n,\|\cdot\|_D)$ embeds in $L_{-k}.$ In this article we use
embeddings in $L_{-p}$ only to state some results in continuous form;  for more 
applications of this concept, see \cite[Ch. 6]{K-book}.

Embeddings in $L_{-p}$ and $k$-intersection bodies admit a Fourier analytic
characterization  that we are going to use throughout
this article:
\begin{pr} \label{posdef} (\cite{K8}, \cite[Th. 6.16]{K-book}) Let $D$ be an origin symmetric star body in $\R^n,\ 0<p<n.$ The space
$(\R^n,\|\cdot\|_D)$ embeds in $L_{-p}$ if and only if the function $\|x\|_D^{-p}$ represents 
a positive definite distribution on $\R^n.$ In particular, $D$ is a $k$-intersection body if
and only if $\|x\|_D^{-k}$ is a positive definite distribution on $\R^n.$
\end{pr}

It was proved in \cite{K6} (see also \cite[Corollary 4.9]{K-book}) that every $n$-dimensional normed space embeds in $L_{-p}$ for
each $p\in [n-3,n).$ In particular, every origin symmetric convex body in $\R^n$ is
a $k$-intersection body for $k=n-3, n-2, n-1.$ On the other hand, the spaces $\ell_q^n,\ q>2$ do not
embed in $L_{-p}$ if $0< p < n-3,$ hence, the unit balls of these spaces are not 
$k$-intersection bodies if $k<n-3$; see \cite{K3}, \cite[Theorem 4.13]{K-book}.  We are going to use a generalization of 
the latter result, the so-called second derivative test for $k$-intersection bodies and 
embeddings in $L_{-p}$, which was first proved for intersection bodies in \cite{K4}
and then generalized in \cite[Theorems 4.19, 4.21]{K-book}.  Recall that for normed 
spaces $X$ and $Y$ and 
$q\in \R,\ q\ge 1,$ the $q$-sum $(X\oplus Y)_q$ of $X$ 
and $Y$ \index{$(X\oplus Y)_q$, q-sum of normed spaces} \index{q-sum of normed spaces}
is defined as the space of pairs $\{(x,y):\ x\in X, y\in Y\}$ with the norm
$$\|(x,y)\|= \left(\|x\|_X^q+\|y\|_Y^q \right)^{1/q}.$$ 

\begin{pr} \label{sdt} Let $n\ge 3,\ k\in {\mathbb N}\cup \{0\}, \ q>2$ and 
let $Y$ be a finite dimensional normed
space of dimension greater or equal to $n.$ Then the $q$-sum of $\R$ and $Y$ 
does not embed in $L_{-p}$ with $0<p<n-2.$ In particular, this direct sum is not
a $k$-intersection body for any $1\le k < n-2.$
\end{pr}

\bigbreak
Let $1\le k<n$ and let $H$ be an $(n-k)$-dimensional subspace of $\R^n.$
Fix any orthonormal basis $e_1,...,e_k$  in the orthogonal subspace $H^\bot.$ 
For a convex body $D$ in $\R^n,$ define the $(n-k)$-dimensional parallel section 
function $A_{D,H}$ as a  function on $\R^k$ such that
$$A_{D,H}(u) = \vol_{n-k}(D\cap \{H+ u_1e_1+...+u_ke_k\})$$
 \begin{equation} \label{paral}
 = \int_{\{x\in \R^n: (x,e_1)=u_1,...,
(x,e_k)=u_k\}} \chi(\|x\|_D)\ dx, \quad u\in \R^k.
\end{equation}

Let $|\cdot|_2$ be the Euclidean norm on $\R^k.$ For every 
$q\in \C,$ the value of the distribution 
$|u|_2^{-q-k}/\Gamma(-q/2)$ on a test function $\phi\in {\cal S}(\R^k)$
can be defined in the usual way (see \cite[p.71]{GS}) and represents an
entire function of $q\in \C.$ If $D$ is infinitely smooth, the
function $A_{D,H}$ is infinitely differentiable at the origin (see \cite[Lemma 2.4]{K-book}), and 
the same regularization procedure can be applied to define the 
action of these distributions on the function $A_{D,H}.$ The function
\begin{equation} \label{fraclapl}
q\mapsto \Big\langle {{|u|_2^{-q-k}}\over {\Gamma(-q/2)}},
A_{D,H}(u) \Big\rangle
\end{equation}
is an entire function of $q\in \C.$ In particular, if $q<0$ 
$$\Big\langle \frac{|u|_2^{-q-k}}{\Gamma(-q/2)}, A_{D,H}(u) \Big\rangle = \frac{1}{\Gamma(-q/2)}
\int_{\R^k} |u|_2^{-q-k} A_{D,H}(u)\ du.$$
If  $q=2m,\ m\in \N\cup \{0\},$ then
$$ \Big\langle {{|u|_2^{-q-k}}\over {\Gamma(-q/2)}}\Big\vert_{q=2m},
A_{D,H}(u) \Big\rangle$$
\begin{equation} \label{intaction}
 = \frac{(-1)^m |S^{k-1}|}{2^{m+1} k(k+2)...(k+2m-2)} \Delta^m A_{D,H}(0),
\end{equation}
where $|S^{k-1}| = 2\pi^{k/2}/\Gamma(k/2)$ is the surface area of the unit
sphere $S^{k-1}$ in $\R^k,$ and 
$\Delta = \sum_{i=1}^k \partial^2/\partial u_i^2$ is
the $k$-dimensional Laplace operator (for details, see \cite[p.71-74]{GS}). 
Since the body $D$ is origin-symmetric, the function $A_{D,H}$ is even,
and for $0<q<2$ we have (see also \cite[p. 49]{K-book}) 
$$\Big\langle \frac{|u|_2^{-q-k}}{\Gamma(-q/2)}, A_{D,H}(u) \Big\rangle $$
\begin{equation} \label{action}
= \frac{1}{\Gamma(-q/2)}\int_{S^{n-1}} \left( \int_0^\infty \frac{A_{D,H}(t\theta) - A_{D,H}(0)}{t^{1+q}}\ dt \right)\ d\theta.
\end{equation}
Note that the function (\ref{fraclapl}) is equal (up to a constant) to the fractional power
of the Laplacian $\Delta^{q/2}A_{D,H}.$

The following proposition was proved in \cite[Th. 2]{K8}. We reproduce the proof here for the 
sake of completeness. We use  a well-known formula (see for example \cite[p. 76]{GS}):
for any $v\in \R^k$ and $q< -k+1,$
$$(v_1^2+...+v_k^2)^{(-q-k)/2} $$
\begin{equation} \label{sph}
= {{\Gamma(-q/2)}\over {2\Gamma((-q-k+1)/2)
\pi^{(k-1)/2}}} \int_{S^{k-1}} |(v,u)|^{-q-k}\ du.
\end{equation}
\begin{pr} \label{main} Let $D$ be an infinitely smooth origin symmetric
convex body in $\R^n$ and $1\le k <n.$ Then for every $(n-k)$-dimensional
subspace $H$ of $\R^n$ and any $q\in \R,\ -k<q< n-k,$
$$\Big\langle {{|u|_2^{-q-k}}\over {\Gamma(-q/2)}},
A_{D,H}(u) \Big\rangle$$
\begin{equation} \label{form1}
 = {{2^{-q-k}\pi^{-k/2}}\over {\Gamma((q+k)/2)(n-q-k)}}
\int_{S^{n-1}\cap H^\bot} \big(\|x\|_D^{-n+q+k}\big)^\wedge(\theta)
\ d\theta.
\end{equation}
Also for every $m\in {\mathbb{N}} \cup \{0\},\ m<(n-k)/2,$
\begin{equation} \label{form2}
 \Delta^m A_{D,H}(0) = {{(-1)^m}\over {2^k\pi^k(n-2m-k)}}
\int_{S^{n-1}\cap H^\bot} (\|x\|_D^{-n+2m+k})^\wedge(\eta)\ d\eta,
\end{equation}
where, as before,  $\Delta$ is the Laplacian on $\R^k.$
\end{pr}

\pf   Let first $q\in (-k,-k+1).$ Then,
$$ \Big\langle {{|u|_2^{-q-k}}\over {\Gamma(-q/2)}},
A_{D,H}(u) \Big\rangle = {1\over {\Gamma(-q/2)}}
\int_{\R^k} |u|_2^{-q-k} A_{D,H}(u)\ du.$$
Using the expression (\ref{paral}) for the function $A_{D,H}$, writing 
the integral in polar coordinates and then using (\ref{sph}), we
see that the right-hand side of the latter equation is equal to 
$${1\over {\Gamma(\frac{-q}{2})}} \int_{\R^n} \big((x,e_1)^2+...
+(x,e_k)^2)^{(-q-k)/2} \chi(\|x\|_D)\ dx=$$
$${1\over {\Gamma(\frac{-q}{2})(n-q-k)}} \int_{S^{n-1}} \left((\theta,e_1)^2+...
+(\theta,e_k)^2 \right)^{(-q-k)/2} \|\theta\|_D^{-n+q+k}\ d\theta=$$
$${1\over {2\Gamma(\frac{-q-k+1}{2})\pi^{\frac{k-1}{2}}(n-q-k)}}\ \times$$
$$\int_{S^{n-1}} \|\theta\|_D^{-n+q+k} \left(
\int_{S^{k-1}} \big| (\sum_{i=1}^k u_ie_i, \theta) \big|^{-q-k}\ du \right)\ d\theta=$$
$${1\over {2\Gamma(\frac{-q-k+1}{2})\pi^{\frac{k-1}{2}}(n-q-k)}}\times$$
\begin{equation} \label{comput}
\int_{S^{k-1}} \left(
\int_{S^{n-1}} \|\theta\|_D^{-n+q+k} \big| (\sum_{i=1}^k u_ie_i, \theta) 
\big|^{-q-k}\ d\theta\right) du.
\end{equation}

Let us show that the function under the integral over $S^{k-1}$
is the Fourier transform of $\|x\|_D^{-n+q+k}$ at the point 
$\sum u_ie_i$. For any even test function $\phi\in {\cal S}(\R^n),$  
using the well-known connection between the Fourier and Radon transforms 
(see \cite[p. 27]{K-book}) and the
expression for the Fourier transform of the distribution 
$|z|^{q+k-1}$ (see \cite[p. 38]{K-book}), we get

$$\langle (\|x\|_D^{-n+q+k})^\wedge, \phi \rangle =
\langle \|x\|_D^{-n+q+k}, \hat\phi \rangle =
\int_{\R^n} \|x\|_D^{-n+q+k} \hat\phi(x)\ dx =$$
$$\int_{S^{n-1}} \|\theta\|_D^{-n+q+k}\left(
\int_0^\infty z^{q+k-1} \hat\phi(z\theta)\ dz\right) d\theta=$$
$${1\over 2} \int_{S^{n-1}} \|\theta\|_D^{-n+q+k}
\Big\langle |z|^{q+k-1}, \hat\phi(z\theta) \Big\rangle \ d\theta=$$
$${{2^{q+k}\sqrt{\pi}\ \Gamma((q+k)/2)}\over 
{2\Gamma((-q-k+1)/2)}} 
\int_{S^{n-1}} \|\theta\|_D^{-n+q+k}
\Big\langle |t|^{-q-k}, 
\int_{(y,\theta)=t} \phi(y)\ dy \Big\rangle \ d\theta=$$
$$ {{2^{q+k}\sqrt{\pi} \Gamma((q+k)/2)}\over 
{2\Gamma((-q-k+1)/2)}} \int_{\R^n} \Big(\int_{S^{n-1}}
|(\theta, y)|^{-q-k} \|\theta\|_D^{-n+q+k}\ d\theta \Big) \phi(y)\ dy.$$
Since $\phi$ is an arbitrary test function, this 
proves that, for every $y\in \R^n\setminus \{0\},$
$$\big(\|x\|_D^{-n+q+k}\big)^\wedge(y) =
{{2^{q+k}\sqrt{\pi} \Gamma((q+k)/2)}\over 
{2\Gamma((-q-k+1)/2)}} \int_{S^{n-1}}
|(\theta, y)|^{-q-k} \|\theta\|_D^{-n+q+k}\ d\theta.$$
Together with (\ref{comput}), the latter equality shows that 
\begin{equation} \label{fracform}
\Big\langle {{|u|_2^{-q-k}}\over {\Gamma(-q/2)}},
A_{D,H}(u) \Big\rangle
\end{equation}
$$ = {{2^{-q-k}\pi^{-k/2}}\over {\Gamma((q+k)/2)(n-q-k)}}
\int_{S^{n-1}\cap H^\bot} \big(\|x\|_D^{-n+q+k}\big)^\wedge(\theta)
\ d\theta,$$
because in our notation $S^{k-1}=S^{n-1}\cap H^\bot.$

We have proved (\ref{fracform}) under the assumption that $q\in (-k,-k+1).$ 
However, both sides of (\ref{fracform})
are analytic functions of $q\in \C$ in the domain where 
$-k<Re(q)<n-k.$   
This implies that the equality (\ref{fracform}) holds for every $q$ from this domain
(see \cite[p. 61]{K-book} for the details of a similar argument).

Putting $q=2m,\ m\in \N\cup \{0\},\ m< (n-k)/2$ in (\ref{fracform}) and applying
(\ref{intaction}) and the fact that $\Gamma(x+1)=x\Gamma(x)$, we get the second
formula.
\endpf
\medbreak

Brunn's theorem
(see for example \cite[Th. 2.3]{K-book}) states that the central hyperplane
section of an origin symmetric convex body has maximal $(n-1)$-dimensional
volume among all hyperplane sections perpendicular to a given direction.
This implies the following
\begin{lemma} \label{laplace} If $D$ is a $2$-smooth origin symmetric convex body in $\R^n,$
then the function $A_{D,H}$ is twice differentiable at the origin and 
$$\Delta A_{D,H}(0) \le 0.$$ Besides that for any $q\in (0,2),$ 
$$\Big\langle \frac{|u|_2^{-q-k}}{\Gamma(-q/2)}, A_{D,H}(u) \Big\rangle \ge 0.$$
\end{lemma}

\pf Differentiability follows from \cite[Lemma 2.4]{K-book}. 
Applying Brunn's theorem to the bodies $D\cap {\rm span}(H,\theta),\ \theta\in S^{n-1}\cap H^\bot$, we see that
the function $t\mapsto A_{D,H}(t\theta)$ has maximum at zero. Therefore,
the interior integral in (\ref{action}) is negative, but $\Gamma(-q/2)<0$ for $q\in (0,2),$
which implies the second statement. The first inequality also follows from the fact 
that each of the functions $t\mapsto A_{D,H}(t e_j),\ j=1,...,k$ has maximum at the origin.
\endpf
\medbreak
We often use Lemma 4.10 from \cite{K-book} for the purpose of approximation
by infinitely smooth bodies. For convenience, let us formulate this lemma:
\begin{lemma} \label{smooth} (\cite[Lemma 4.10]{K-book}) \label{intbody-approx} Let $1\le k < n.$
Suppose that $D$ is an origin-symmetric convex body in $\R^n$
that is not a $k$-intersection body. Then there exists a
sequence $D_m$ of origin-symmetric convex bodies so that $D_m$
converges to $D$ in the radial metric, each $D_m$ is infinitely
smooth, has strictly positive curvature and each $D_m$ is not 
a $k$-intersection body.
\end{lemma}
If in addition $D$ is invariant with respect to $R_\theta,$ one can choose $D_m$
with the same property. 

\section{Connection with intersection bodies}

We now return to the complex case. The following simple observation is crucial for
applications of the Fourier methods to convex bodies in the complex case:

\begin{lemma} \label{const} Suppose that $K$ is an origin-symmetric infinitely smooth 
invariant with respect to all  $R_\theta$ star body in $\R^{2n}.$ Then for 
every $0<p<2n$ and $\xi\in S^{2n-1}$ the Fourier transform of the distribution $\|x\|_K^{-p}$ is 
a constant function on $S^{2n-1}\cap H_\xi^\bot.$
\end{lemma}
\pf By \cite[Lemma 3.16]{K-book}, the Fourier transform of $\|x\|_K^{-p}$ is a continuous
function outside of the origin in $\R^{2n}.$ The function $\|x\|_K$ is invariant with respect 
to all $R_\theta$, so by the 
connection between the Fourier transform of distributions and linear transformations,
the Fourier transform of $\|x\|_K^{-p}$ is also invariant with respect to all $R_\theta.$
Recall that the two-dimensional space $H_\xi^\bot$ is spanned by  vectors
$\xi$ and $\xi^\bot$ (see the Introduction). Every vector in $S^{2n-1}\cap H_\xi^\bot$
is the image of $\xi$ under one of the coordinate-wise rotations $R_\theta$, so 
the Fourier transform of $\|x\|_K^{-p}$ is a constant function on $S^{2n-1}\cap H_\xi^\bot.$
\endpf
Of course, this argument also applies to the Fourier transform of any distribution
of the form $h(\|x\|_K).$

Similarly to the real case (see \cite{K1}, \cite[Theorem 3.8]{K-book}), one can express the volume 
of hyperplane sections in terms of the Fourier transform.

\begin{theorem} \label{volume-ft} Let $K$ be an infinitely smooth origin symmetric 
invariant with respect to $R_\theta$ convex 
body in $\R^{2n}, n\ge 2.$  For every $\xi\in S^{2n-1},$ we have
$$\vol_{2n-2}(K\cap H_\xi) = \frac{1}{4\pi(n-1)} \left(\|x\|_K^{-2n+2}\right)^\wedge(\xi).$$
\end{theorem}
 
\pf Let us fix $\xi\in S^{2n-1}.$ We apply formula (\ref{form2}) with $2n$ in place of $n$, 
$H=H_\xi,\ k=2,\ m=0.$ We get 
$$\vol_{2n-2}(K\cap H_\xi) = A_{K,H_\xi}(0)= \frac{1}{8\pi^2(n-1)} \int_{S^{2n-1}\cap H_\xi^\bot}
\left(\|x\|_K^{-2n+2}\right)^\wedge(\eta)\ d\eta.$$
By Lemma \ref{const}, the function under the integral in the right hand side is constant
on the circle $S^{2n-1}\cap H_\xi^\bot.$ Since $\xi\in H_\xi^\bot,$  the integral is equal to 
$2\pi \left(\|x\|_K^{-2n+2}\right)^\wedge(\xi).$
\endpf

The connection between the complex Busemann-Petty problem
and intersection bodies is as follows: 
\begin{theorem} \label{connection} The answer to the complex Busemann-Petty problem in $\C^n$
is affirmative if and only if every origin symmetric invariant with respect to all $R_\theta$ 
convex body in $\R^{2n}$ is a $2$-intersection body.
\end{theorem}

This theorem will follow from the next two lemmas. Note that, since we can approximate
the body $K$ in the radial metric from inside by infinitely smooth convex bodies invariant 
with respect to all $R_\theta,$ and also approximate $L$ from outside in the same way,
we can argue that
if the answer to the complex Busemann-Petty problem is affirmative 
for infinitely smooth bodies $K$ and $L$ then it is affirmative in general.

\begin{lemma} \label{con1} Let $K$ and $L$ be infinitely smooth invariant with respect to $R_\theta$ 
convex bodies in $\R^{2n}$ so that 
$K$ is a $2$-intersection body and, for every $\xi\in S^{2n-1},$
$$\vol_{2n-2}(K\cap H_\xi) \le \vol_{2n-2}(L\cap H_\xi).$$
Then
$$\vol_{2n}(K) \le \vol_{2n}(L).$$ 
\end{lemma}

\pf By \cite[Lemma  3.16]{K-book}, the Fourier transforms of the distributions $\|x\|_K^{-2n+2}$,
$\|x\|_L^{-2n+2}$ and $\|x\|_K^{-2}$ are continuous functions outside of the origin in $\R^{2n}.$
By Theorem \ref{volume-ft} and Proposition \ref{posdef}, the conditions of the lemma imply that
for every $\xi\in S^{2n-1},$
$$ \left(\|x\|_K^{-2n+2}\right)^\wedge(\xi) \le \left(\|x\|_L^{-2n+2}\right)^\wedge(\xi)$$
and
$$\left(\|x\|_K^{-2}\right)^\wedge(\xi)\ge 0.$$
Therefore,
$$\int_{S^{2n-1}} \left(\|x\|_K^{-2n+2}\right)^\wedge(\xi) \left(\|x\|_K^{-2}\right)^\wedge(\xi)\  d\xi$$
$$ \le \int_{S^{2n-1}} \left(\|x\|_L^{-2n+2}\right)^\wedge(\xi) \left(\|x\|_K^{-2}\right)^\wedge(\xi)\  d\xi.$$
Now we apply Parseval's formula on the sphere, Proposition \ref{parseval}, to
remove the Fourier transforms in the latter inequality and then use the polar formula for the volume
and H\"older's inequality:
$$2n\ \vol_{2n}(K) = \int_{S^{2n-1}} \|x\|_K^{-2n} \  dx
\le \int_{S^{2n-1}} \|x\|_L^{-2n+2}  \|x\|_K^{-2}\  dx$$
$$ \le  \left(\int_{S^{2n-1}} \|x\|_L^{-2n} \  dx\right)^{\frac{n-1}{n}} 
\left(\int_{S^{2n-1}} \|x\|_K^{-2n} \  dx\right)^{\frac{1}{n}}$$
$$ = \left(2n\ \vol_{2n}(L)\right)^{\frac{n-1}{n}}
 \left(2n\ \vol_{2n}(K)\right)^{\frac{1}{n}},$$
which gives the result.
\endpf

\begin{lemma} \label{con2} Suppose that there exists an origin symmetric invariant with
respect to all $R_\theta$ convex body $L$ in $\R^{2n}$
which is not a $2$-intersection body. Then one can perturb $L$ twice to construct other origin symmetric
invariant with respect to $R_\theta$ convex bodies $L^{'}$ and $K$ in $\R^{2n}$ such that for every $\xi\in S^{2n-1},$
$$\vol_{2n-2}(K\cap H_\xi) \le \vol_{2n-2}(L^{'} \cap H_\xi),$$
but
$$\vol_{2n}(K) > \vol_{2n}(L^{'}).$$ 
\end{lemma}

\pf We can assume that the body $L$ is infinitely smooth and has strictly positive curvature. In fact,
approximating $L$ in the radial metric by infinitely smooth invariant with respect to all $R_\theta$ 
convex bodies with strictly positive curvature, we 
get by Lemma \ref{smooth} that approximating bodies can not all be $2$-intersection
bodies. So there exists an infinitely smooth invariant with respect to all $R_\theta$ convex body  $L^{'}$
with strictly positive curvature that is not 
a $2$-intersection body. 

Now as $L$ is infinitely smooth, by \cite[Lemma 3.16]{K-book}, the Fourier transform of $\|x\|_L^{-2}$
is a continuous function outside of the origin in $\R^{2n}.$ The body $L$ is not a $2$-intersection body,
so by Proposition \ref{posdef}, the Fourier transform $\left(\|x\|_L^{-2}\right)^\wedge$  is negative on some 
open subset $\Omega$ of the sphere $S^{2n-1}.$ 

Since $L$ is invariant with respect to rotations $R_\theta,$
we can assume that the set $\Omega$ is also invariant with respect to rotations $R_\theta.$ This allows us 
to choose an even non-negative invariant with respect to rotations $R_\theta$ function $f\in C^\infty(S^{2n-1})$
which is supported in $\Omega.$ Extend $f$ to an even homogeneous function $f(x/|x|_2)|x|_2^{-2}$ of degree -2
on $\R^{2n}.$ By \cite[Lemma 3.16]{K-book}, the Fourier transform of this extension is an even homogeneous 
function of degree -2n+2 on $\R^{2n},$ whose restriction to the sphere is infinitely smooth: 
$$\left(f(x/|x|_2)|x|_2^{-2}\right)^\wedge(y) = g(y/|y|_2)|y|_2^{-2n+2},$$
 where $g\in C^\infty(S^{2n-1}).$ By the connection
between the Fourier transform and linear transformations, the function $g$ is also invariant with respect
to rotations $R_\theta.$

Define a body $K$ in $\R^{2n}$ by
\begin{equation} \label{counter}
\|x\|_K^{-2n+2} = \|x\|_L^{-2n+2} - \epsilon g(x/|x|_2)|x|_2^{-2n+2}.
\end{equation}
For small enough $\epsilon$ the body $K$ is convex. This essentially follows from a simple 
two-dimensional argument: if $h$ is a strictly concave function on an interval $[a,b]$ and
$u$ is a twice differentiable function on $[a,b]$, then for small $\epsilon$ the function
$h+\epsilon u$ is also concave. Note that here we use the condition that $L$ has strictly
positive curvature. Besides that, the body $K$ is invariant with respect to rotations $R_\theta$
because so are the body $L$ and the function $g.$ We can now choose $\epsilon$ so that 
$K$ is an origin symmetric  invariant with respect to all $R_\theta$ convex body in $\R^{2n}.$

Let us prove that the bodies $K$ and $L$ provide the necessary counterexample.
We apply the Fourier transform to both sides of (\ref{counter}).  By 
definition of the function $g$ and since $f$ is non-negative,
we get that for every $\xi\in S^{2n-1}$
$$\left(\|x\|_K^{-2n+2}\right)^\wedge(\xi) =\left(\|x\|_L^{-2n+2}\right)^\wedge(\xi) - (2\pi)^{2n}\epsilon f(\xi)
\le \left(\|x\|_L^{-2n+2}\right)^\wedge(\xi).$$
By Theorem \ref{volume-ft}, this means that for every $\xi$
$$\vol_{2n-2}(K\cap H_\xi) \le \vol_{2n-2}(L\cap H_\xi).$$
On the other hand, the function $f$ is positive only where $\left(\|x\|_L^{-2}\right)^\wedge$
is negative, so 
$$\int_{S^{2n-1}} \left(\|x\|_K^{-2n+2}\right)^\wedge(\xi) \left(\|x\|_L^{-2}\right)^\wedge(\xi)\  d\xi$$
$$ = \int_{S^{2n-1}} \left(\|x\|_L^{-2n+2}\right)^\wedge(\xi) \left(\|x\|_L^{-2}\right)^\wedge(\xi)\  d\xi$$
$$- (2\pi)^{2n}\epsilon \int_{S^{2n-1}} \left(\|x\|_L^{-2}\right)^\wedge(\xi) f(\xi)\ d\xi$$
$$> \int_{S^{2n-1}} \left(\|x\|_L^{-2n+2}\right)^\wedge(\xi) \left(\|x\|_L^{-2}\right)^\wedge(\xi)\  d\xi.$$
The end of the proof is similar to that of the previous lemma - we apply Parseval's 
formula to remove Fourier transforms and then use H\"older's inequality and the polar 
formula for the volume to get $\vol_n(K)>\vol_n(L).$
\endpf

\section{The solution of the problem}

It is known (see \cite{K6} or \cite[Corollary 4.9]{K-book} plus Proposition \ref{posdef}) that for every origin symmetric convex 
body $K$ in $\R^{2n}, n\ge 2$ the space $(\R^{2n},\|\cdot\|_K)$ embeds in $L_{-p}$ for each $p\in [2n-3,2n),$ 
or, in other words,  every origin-symmetric 
convex body  in $\R^{2n}$ is a $(2n-3)$-, $(2n-2)$- and $(2n-1)$-intersection body. On the other hand, 
for $q>2$ the unit ball of the real space $\ell_q^{2n}$ is not a $(2n-4)$-intersection body ,
and, moreover,  $\R^{2n}$ provided with the norm of this space does not embed in $L_{-p}$ with $p<2n-3$
(see \cite{K3} or  \cite[Th. 4.13]{K-book}).

Now we have to find out what happens if we consider convex bodies invariant with respect to all $R_\theta.$
It immediately follows from the second derivative test (\cite[Th. 4.19 and 4.21]{K-book} ;
see Corollary \ref{bqn} below) that for $q>2$ the complex space $\ell_q^n$ does not embed in $L_{-p}$ with
$p<2n-4,$ which means that the unit ball $B_q^n$ of this space (which is invariant with respect
to all $R_\theta)$ is not a $k$-intersection 
body with $k<2n-4.$ The only  question that remains open is what happens in the interval $p\in [2n-4,2n-3).$ 
The following result answers this question.

\begin{theorem} \label{2n-4}Let  $n\ge 3.$ Every origin symmetric invariant with respect to $R_\theta$ 
convex body $K$ in $\R^{2n}$ is a $(2n-4)$-intersection body. Moreover, the space $(\R^{2n},\|\cdot\|_K)$
embeds in $L_{-p}$ for every $p\in [2n-4,2n).$

If $n=2$ the space $(\R^{2n},\|\cdot\|_K)$ embeds in $L_{-p}$ for every $p\in (0,4).$
\end{theorem}

\pf  By Lemma \ref{smooth}, it is enough to prove the result in the case where $K$ is infinitely 
smooth. Fix $\xi\in S^{2n-1}.$ 

Let $n\ge 3.$ Applying formula (\ref{form2}) and then Lemma \ref{const} with 
$ H=H_\xi,\ m=1,\ k=2$ and dimension $2n$ instead of $n,$ we get
$$\Delta A_{K,H_\xi}(0) = \frac{-1}{8\pi^2(n-2)} \int_{S^{n-1}\cap H_\xi^\bot} 
(\|x\|_K^{-2n+4})^\wedge(\eta)\ d\eta$$
$$= \frac{-2\pi}{8\pi^2(n-2)} \left(\|x\|_K^{-2n+4}\right)^\wedge(\xi).$$
By Brunn's theorem (see Lemma \ref{laplace}),  $\left(\|x\|_K^{-2n+4}\right)^\wedge(\xi)\ge 0$  
for every $\xi\in S^{2n-1},$ so 
$\|x\|_K^{-2n+4}$ is a positive definite distribution on $\R^{2n}.$ By Proposition \ref{posdef}, 
$K$ is a $(2n-4)$-intersection body.
 
Now let $n\ge 2.$ For $0<q<2$, formula (\ref{form1}) and Lemma \ref{laplace} imply
that  $\left(\|x\|_K^{-2n+q+2}\right)^\wedge(\xi)\ge 0.$ By Proposition \ref{posdef}, the space
$(\R^{2n},\|\cdot\|_K)$ embeds in $L_{-2n+q+2},$ and, using the range of $q$, every such 
space embeds in $L_{-p},\ p\in (2n-4,2n-2).$ As mentioned before, these spaces also embed
in $L_{-p},\ p\in [2n-3,2n),$ because so does any $2n$-dimensional normed space.
 \endpf

We now give an example of an origin symmetric invariant with respect to all $R_\theta$ 
convex body in $\R^{2n}$ which is not a $k$-intersection body for any $1\le k < 2n-4.$

Denote by $B_n^q$ the unit ball of the complex space $\ell_q^n$ considered as 
a subset of $\R^{2n}:$
$$B_q^n = \{\xi\in \R^{2n}: \ \|\xi\|_q=\left(\left(\xi_{11}^2+\xi_{12}^2\right)^{q/2}+...
+ \left(\xi_{n1}^2+\xi_{n2}^2\right)^{q/2}\right)^{1/q}\le 1\}.$$
If $q\ge 1$ then $B_q^n$ is an origin symmetric invariant  with respect to $R_\theta$
convex body in $\R^{2n}.$

The next theorem immediately follows from Proposition \ref{sdt}.

\begin{theorem} \label{bqn} If $q>2$ then the space $(\R^{2n}, \|\cdot\|_q)$ does not embed in $L_{-p}$
with $0<p<2n-4.$ In particular, the body $B_q^n$ is not a $k$-intersection body for
any $1\le k < 2n-4.$
\end{theorem}

\pf The space $(\R^{2n},\|\cdot\|_q)$ contains as a subspace the $q$-sum of $\R$
and a $(2n-2)$-dimensional subspace $(\R^{2n-2},\|\cdot\|_q).$ This $q$-sum 
does not embed in $L_{-p},\ 0<p<2n-4$ by Proposition \ref{sdt}. By a result of
E.Milman \cite{Mi}, the larger space cannot embed in  $L_{-p},\ 0<p<2n-4$ either
(the proof in \cite{Mi} is only for integers $p$, but it is exactly the same for
non-integers; note that for the complex Busemann-Petty problem 
we need only the second statement of the corollary, where $p$ is an integer).
\endpf

We are now ready to prove the main result of this article:

\begin{theorem} The solution to the complex Busemann-Petty problem in $\C^n$ 
is affirmative if $n\le 3$ and it is negative if $n\ge 4.$
\end{theorem}

\pf By Theorem \ref{2n-4}, every origin symmetric invariant with respect to $R_\theta$ convex
body in $\R^6$ (where $n=3$) is a $2n-4=2$-intersection body, and in $\R^4$ (where $n=2$) it 
is a $2n-2=2$-intersection body.
The affirmative answers for $n=3$ and $n=2$ follow now from Theorem \ref{connection}.

If $n\ge 4$ then $2n-4>2,$ so by Theorem \ref{bqn} the body $B_q^n$ is not a $2$-intersection
body. The negative answer follows from Theorem \ref{connection}.
\endpf

\bigbreak

{\bf Remark 1.} The transition between the dimensions $n=3$ and $n=4$ is due to the fact that 
convexity controls only derivatives of the second order.  To see this let us look again at formula
(\ref{form2}), which we apply with $k=2.$ We want to get information about the Fourier transform 
of $\|x\|_D^{-2},$ so we need to choose $m$ so that $-2n+2m+2=-2.$ If $n=3$ then $m=1,$ but when 
$n=4$ we need $m=2.$ This means that for $n=3$ we consider $\Delta A_{K,H}(0),$ which is always 
negative by convexity, but when $n=4$ we look at $\Delta^2A_{K,H}(0),$ which is not controlled
by convexity and can be sign-changing. One can construct a counterexample in dimension $n=4$ 
using this argument, similarly to how it was done for the ``real" Busemann-Petty problem; 
see \cite[Corollary 4.4]{K-book}.

\medbreak

{\bf Remark 2.} Applying Theorem \ref{2n-4} to $n=2$ we get that every two-dimensional complex
normed space (which is a 4-dimensional real normed space) embeds in $L_{-p}$ for every $p\in [-1,0).$
By  \cite[Th. 6.4]{KKYY}, this implies that every such space embeds isometrically in $L_0.$
The concept of embedding in $L_0$ was introduced in \cite{KKYY}: a normed space
$(\R^n,\|\cdot\|)$ embeds in $L_0$ if there exist a probability measure $\mu$ on $S^{n-1}$
and a constant $C$ so that for every $x\in \R^n,\ x\neq 0$
$$\log \|x\| = \int_{S^{n-1}} \log |(x,\xi)|\ d\mu(\xi)  + C.$$
We have 
\begin{theorem} Every two-dimensional complex normed space embeds in $L_0.$
On the other hand, there exist two-dimensional complex normed spaces that do not
embed isometrically in any $L_p,\ p>0.$
\end{theorem}
An example supporting the second claim is the complex space $\ell_q^2$ with $q>2.$ This 
follows from a version of the second derivative test proved in \cite{KL} (see also 
\cite[Theorem 6.11]{K-book}). Recall that every two-dimensional
real normed space embeds isometrically in $L_1$ (see \cite{Fe}, \cite{He}, \cite{Li} or \cite[p. 120]{K-book}),
but the real space $\ell_q^2$ does not embed isometrically  in any $L_p,\ 1<p\le 2,$ as proved by Dor \cite{Do}; see also \cite[p. 124]{K-book}.

\bigbreak

{\bf Acknowledgments:} 
The first named author was supported in part by the NSF grants DMS-0455696
and DMS-0652571. The third named author was partially supported by the 
NSF grant DMS-0652571 and by the European
Network PHD, FP6 Marie Curie Actions, Contract MCRN-511953.

\end{document}